\documentclass[letterpaper, 10 pt, conference]{ieeeconf} 
\IEEEoverridecommandlockouts 

\usepackage{graphicx}
\usepackage{booktabs}
\usepackage{amsmath} 
\let\labelindent\relax
\usepackage{enumitem}
\usepackage{amsfonts}
\usepackage{amssymb}  
\usepackage{mathtools} 
\usepackage{mathrsfs}
\usepackage{balance}
\usepackage{url}
\usepackage{hyperref}
\usepackage{algorithm}
\usepackage{algpseudocode}
\DeclareMathAlphabet{\mathpzc}{OT1}{pzc}{m}{it}
\usepackage{makecell}
\usepackage[noadjust]{cite}
\usepackage{xcolor}
\usepackage{comment}

\newcommand{\shnote}[1]%
    {\textcolor{magenta}{ #1}}

\newcommand{\zgnote}[1]%
    { \textbf{\textcolor{purple}{#1}} \newline}

\newcommand{\tinit}{t}
\newcommand{\tstep}{\delta \tinit}
\newcommand{\tvar}{s}
\newcommand{\thor}{0} 
\newcommand{\state}{x} 

\newcommand{\tset}{\mathcal{X}} 
\newcommand{\dyn}{f} 
\newcommand{\cdyn}{g} 

\newcommand{\ctrl}{u} 
\newcommand{\csig}{\ctrl(\cdot)} 

\newcommand{\cset}{\mathcal{U}} 
\newcommand{\cfset}{\mathbb{U}} 
\newcommand{\ACSS}{\mathbb{U}_\text{a}} 
\newcommand{\rACSS}{\bar {\mathbb{U}}_\text{a}}
\newcommand{\rACS}{\bar {\mathcal{U}}_\text{a}}

\newcommand{\ACSSo}{\mathbb{U}_{\text{a},1}} 
\newcommand{\ACSSt}{\mathbb{U}_{\text{a},2}} 
\newcommand{\ACSo}{\mathcal{U}_{\text{a},1}} 
\newcommand{\ACSt}{\mathcal{U}_{\text{a},2}} 

\newcommand{\ACS}{\mathcal{U}_\text{a}} 

\newcommand{\traj}{\xi} 
\newcommand{\subt}{\phi} 

\newcommand{\subs}{z} 
\newcommand{\tsset}{\mathcal{Z}} 

\newcommand{\proj}{\mathcal{P}}
\newcommand{\proji}{\mathcal{P}^{-1}}

\newcommand{\valfunc}{V} 
\newcommand{\tclvf}{\valfunc_\gamma} 
\newcommand{\clvf}{\valfunc_\gamma^\infty} 
\newcommand{\rclvf}{\bar{\valfunc}^\infty_\gamma} %
\newcommand{\rtclvf}{\bar{\valfunc}_\gamma} 
\newcommand{\clvfa}{\valfunc_{\gamma ,1}^\infty} 
\newcommand{\clvfb}{\valfunc_{\gamma ,2}^\infty} 

\newtheorem{remark}{Remark}
\newtheorem{proposition}{Proposition}
\newtheorem{definition}{Definition}

\newtheorem{theorem}{Theorem}
\newtheorem{lemma}[theorem]{Lemma}

\title{\LARGE \bf
Synthesizing Control Lyapunov-Value Functions for High-Dimensional Systems Using System Decomposition and Admissible Control Sets
}

\author{Zheng Gong, Hyun Joe Jeong and Sylvia Herbert
\thanks{All authors are with the Department of Mechanical and Aerospace Engineering at the University of California, San Diego 
\{{\href{mailto:zhgong@ucsd.edu}{zhgong}, \href{mailto:hjjeong@ucsd.edu}{hjjeong},\href{mailto:sherbert@ucsd.edu}{sherbert}\}@ucsd.edu.}
}%
}

\begin{document}

\maketitle
\begin{abstract}
Control Lyapunov functions (CLFs) play a vital role in modern control applications, but finding them remains a problem. Recently, the control Lyapunov-value function (CLVF) and robust CLVF  have been proposed as solutions for nonlinear time-invariant systems with bounded control and disturbance. However, the CLVF suffers from the ``curse of dimensionality,'' which hinders its application to practical high-dimensional systems. In this paper, we propose a method to decompose systems of a particular coupled nonlinear structure, in order to solve for the CLVF in each low-dimensional subsystem.  We then reconstruct the full-dimensional CLVF and provide sufficient conditions for when this reconstruction is exact. Moreover, a point-wise optimal controller can be obtained using a quadratic program. We also show that when the exact reconstruction is impossible, the subsystems' CLVFs and their ``admissible control sets'' can be used to generate a Lipschitz continuous CLF. We provide several numerical examples to validate the theory and show computational efficiency. 
\end{abstract}

\section{Introduction}
Ensuring stability is one of the most important tasks for autonomous systems operating in the real world. Control Lyapunov functions (CLFs) are energy-like functions that stabilize systems to their equilibrium points~\cite{sontag1989universal,freeman1996control}. However, it is well known that there lacks a ``universal'' way of finding CLFs for general nonlinear systems, especially with state or input constraints. Hand-designed and task-specific CLFs have been proposed~\cite{clf_zubov,ogren2001control,artstein1983stabilization}. 

A method to construct control Lyapunov-like functions for nonlinear systems with control and disturbance bounds has been proposed \cite{gong2022constructing,gong2024robust}. These control Lyapunov value functions (CLVFs) are based on Hamilton-Jacobi reachability analysis~\cite{evans_hj,bansal2017hamilton}, and are constructed through dynamic programming.
The CLVF guarantees exponential stabilizability, works for general nonlinear dynamics, and deals with state and input constraints well. However, it requires solving a CLVF variational inequality (CLVF-VI) on a discrete grid iteratively, therefore suffering from the ``curse of dimensionality.'' In the HJ reachability community, many works are proposed to solve this issue, including incorporating reinforcement learning~\cite{fisac2019bridging}, deep learning~\cite{bansal2021deepreach}, self-contained subsystems decomposition (SCSD)~\cite{chen2017exact,chen2018decomposition}, warming starting~\cite{herbert2019warm} and model reduction~\cite{liu2020leveraging,holmes2020reachable}. 

The SCSD method provides guaranteed exact solutions for a Hamilton-Jacobi reachability problem under certain assumptions on the coupled nonlinear system and the problem formulation. This approach was further generalized recently through the use of an ``admissible control signal set'' (ACSS) and ``admissible control set'' (ACS)~\cite{zheng2023decomp}, which are used to forward propagate and refine the reconstructed value function. 
This line of work has only been applied to reachability problems, and only provides guarantees on the level sets (reachable sets) of the computed value function, rather than guarantees on the value function itself.

\begin{figure}[t]
    \centering
    \includegraphics[width=\columnwidth]{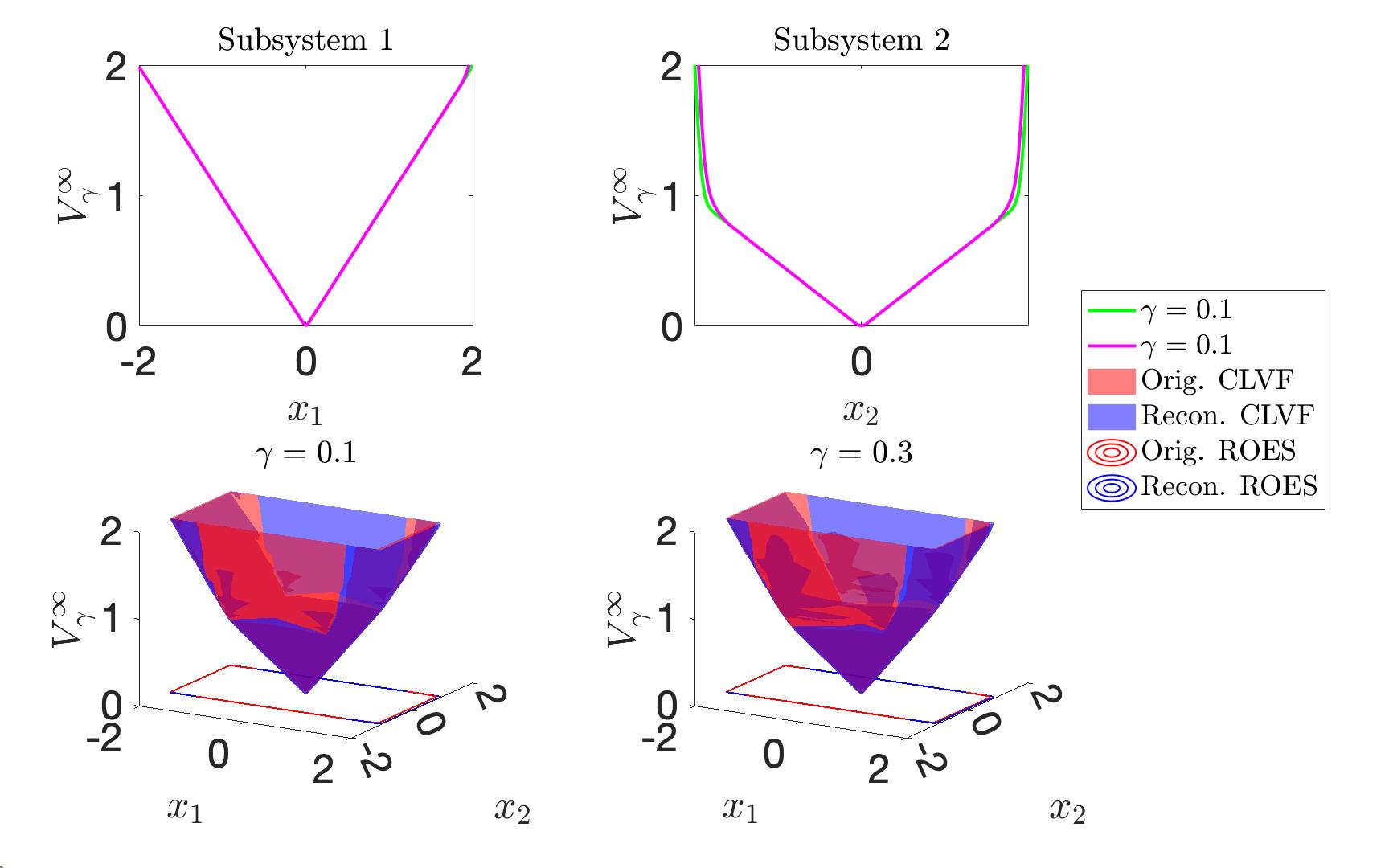}
    \caption{Original and reconstructed CLVFs for a nonlinear 2D system~\eqref{eqn: nonlinear2D}. Top, CLVFs for the two subsystems with two different rates of exponential stabilizability ($\gamma)$. Bottom, comparison of the reconstructed CLVFs and original CLVFs for different $\gamma$, along with their regions of exponential stabilizability (ROES). Both reconstructed CLVFs are identical to the original CLVFs, and the corresponding ROESs also match, validating Lemma~\ref{lemma: exact_decomposition}. }
    \label{fig:2Dsys}
\end{figure}

In this work, we generalize the SDSC and ACS decomposition work to numerically compute CLVFs for high dimensional nonlinear systems with input constraints. The proposed method first applies the standard SCSD and finds the CLVFs for the subsystems. A value function in the original state space is then reconstructed using the subsystems' CLVFs, and the ACS is used to determine the domain where this reconstruction is exact. This is demonstrated in Fig.~\ref{fig:2Dsys}. For systems in which exact reconstruction is not possible, the subsystems' CLVFs and ACSs can be used to generate a Lipschitz continuous CLF and stabilize the original system. The improved scalability of the proposed method is validated in several numerical examples.

The main contributions of this work are:
\begin{enumerate}
    \item We extend the definition of the ACSS and ACS to multi-input systems' CLVFs. 
    \item We propose a method that combines the SCSD, ACS, and CLVF, to reconstruct a value function for high-dimensional systems whose CLF is not trivial to obtain,
    \item We provide sufficient conditions on when the reconstruction is exact, i.e., the reconstructed value function is the CLVF for the original system.
    \item For the case when the exact reconstruction is impossible, we show that a Lipschitz continuous CLF can be reconstructed using the subsystems' CLVFs.
\end{enumerate}

\section{Background}

Consider the following nonlinear control-affine system:
\begin{align}\label{eqn:dyn_sys}
    \frac{d \state}{d \tvar} = \dot \state  = \dyn (\state)+ \cdyn(\state) \cdot \ctrl, \hspace{3mm} \state(0) = \state _0, \hspace{3mm} \tvar \in [\tinit,\thor],
\end{align}
where $\tvar$ is the time, $\state \in \tset \subseteq \mathbb R^n $ is the state, and $\ctrl \in \cset \subseteq \mathbb R^m$ is the control input, where $\cset$ is a compact set. Assume $\dyn: \tset \mapsto \mathbb R^n$ and $\cdyn:\tset \mapsto  \mathbb R^n$ are uniformly continuous, bounded, and Lipschitz continuous in $\state$, and the origin is one equilibrium point, i.e., $\dyn(0)+\cdyn(0)\cdot 0 = 0$. Further assume the control signal $\ctrl (\cdot)$ is drawn from measurable control functions: 
\begin{align*}
    \ctrl(\cdot) \in \cfset := \{ \ctrl: [\tinit,\thor] \mapsto \cset, \ctrl (\cdot) \text{ is measurable}\} .
\end{align*}
Under these assumptions, we can solve for a unique solution of \eqref{eqn:dyn_sys}, denoted as $\traj (\tvar ;\tinit,\state,\ctrl (\cdot))$ (in short $\traj(\tvar)$).

\subsection{CLF and CLVF}
\begin{definition} Given the following time-varying CLVF (TV-CLVF) 
\begin{align} \label{eqn:tv_clvf}
    \tclvf (\state,\tinit) = \min _{\ctrl \in \cfset}  \max_{ \tvar \in [\tinit , \thor]} e^{\gamma( \tvar -\tinit)}  \ell ( \traj (\tvar) ) ,
\end{align}
the CLVF $\clvf: D_\gamma \mapsto \mathbb R $ of \eqref{eqn:dyn_sys} is the limit function: 
\begin{align} \label{eqn:clvf}
    \clvf (\state) = \lim _{\tinit \rightarrow -\infty} \tclvf (\state,\tinit).
\end{align} 
Here, $\gamma$ is a user-specified parameter that represents the desired exponential decay rate. The domain of the CLVF is given as $ D_\gamma \subseteq \mathbb R^n$.
\end{definition}

Proposition 2 of~\cite{gong2022constructing} shows there exists $\ctrl \in \cset$ s.t. 
\begin{align} \label{eqn:Vdot}
    \dot V_\gamma^\infty (\state) \leq \gamma \clvf(\state), \quad \forall \state \in \mathcal D_\gamma. 
\end{align}

Note that since we assume the system has an equilibrium point, the definition of CLVF is simplified compared to \cite{gong2022constructing}. The CLVF value of a state captures the largest exponentially amplified deviation between the origin and the trajectory starting from that state. If this is finite, then there must exist some control signal for which the system trajectory converges to the origin at the same rate $\gamma$. The domain of the CLVF $D_\gamma $ is also known as the region of exponential stabilizability (ROES). A larger $\gamma$ results in a faster convergence to the origin, but a smaller ROES. 

Under the scope of this paper, the CLVF is a Lipschitz continuous CLF, while the main difference is that using CLVF, the user can trade off between a faster convergence or a larger ROES. 

It has been shown that the TV-CLVF is the unique viscosity solution to the following VI, 
\begin{align*} 
    &0 = \max\biggl\{\ell(\state) - \tclvf(\state,
    \tinit), \notag \\
    &\hspace{.5em}  D_\tinit   \tclvf + \min_{ \ctrl \in \cset} D_\state \tclvf \cdot (\dyn (\state)+ \cdyn(\state) \cdot \ctrl )+ \gamma  \tclvf (\state,\tinit) \biggl\},
\end{align*} 
with terminal condition $\tclvf (\state,\thor) = \ell(\state)$. The CLVF can be obtained by solving this VI backward in time until convergence. Here, $ D_\tinit$ denotes the time derivative, and $ D_\state$ denotes the gradient with respect to $\state$.

\subsection{System decomposition and ACS}
\begin{definition}
\label{def: SCSD}
(SCSD) Given system \eqref{eqn:dyn_sys} and assume there exists state partitions $\subs_1 = (\state_1,\state_c) \in \tsset_1$, $\subs_2 = (\state_2,\state_c) \in \tsset_2$, where $\state_1\in\mathbb{R}^{n_1}$, $\state_2\in\mathbb{R}^{n_2}$, $\state_c\in\mathbb{R}^{n_c}$, $n_1,n_2>0$, $n_c\ge 0$, $n_1+n_2+n_c=n$. Assume also the control inputs can be partitioned similarly with $v_1 = (u_1,u_c) \in \mathcal V_1$, $v_2 = (u_2,u_c) \in \mathcal V_2$, where $u_1 \in \mathbb R^{m_1}$, $u_2 \in \mathbb R^{m_2}$, $u_c \in \mathbb R^{m_c}$ and $m_1+m_2+m_c = m$. The two self-contained subsystems of \eqref{eqn:dyn_sys} are
\begin{align} \label{eqn:subsys}
     \dot \subs_1 = \dyn_1(\subs_1)+g_1(\subs_1)v_1, \quad
    \dot \subs_2 = \dyn_2(\subs_2)+g_2(\subs_2)v_2,
\end{align}
with corresponding solution $\subt_1 (\tvar)$, $\subt_2 (\tvar)$.
\end{definition}

When $\state_c$ or/and $\ctrl_c$ exist, we say the subsystems are coupled through the common states $\state_c$ or/and controls $\ctrl_c$, and they have shared states or/and controls. The shared controls are one of the main causes of the inconsistency between the reconstructed value function from the subsystems and the original value function.

\begin{remark} \label{remark:shared_state}
    After decomposition, if the subsystems only have shared states, but no shared controls, then the shared states must have dynamics $\dot \state_c = \dyn_c (\state_c)$. This means the control has no impact on the shared states. In other words, starting from the same initial condition, the trajectory is the same for all subsystems. 
    On the other hand, if the subsystems have shared control, they must also have shared states. 
\end{remark}

We introduce the projection and back projection operator on both state and control. 
\begin{definition}
\label{def:proj} (Projection and back projection operators.)\\
A \textbf{state projection operator} $\proj_{\state, i}: \mathbb R^n \mapsto \mathbb R^{n_i+n_c}$ maps a state $\state \in \mathbb R^n$ to a state $\subs_i \in \mathbb R^{n_i+n_c}$:
\begin{align*} 
    \proj_{\state,i} (\state) := (\state_i,\state_c) = \subs_i.
\end{align*}
A \textbf{control projection operator} $\proj_{\ctrl, i}: \mathbb R^m \mapsto \mathbb R^{m_i+m_c}$ maps a control $\ctrl \in \mathbb R^m$ to a control $v_i \in \mathbb R^{m_i+m_c}$:
\begin{align*} 
    \proj_{\ctrl,i} (\ctrl) := (\ctrl_i,\ctrl_c) = v_i.
\end{align*}
A \textbf{state back projection operator} $\proji_{\state,i}: \mathbb R^{n_i+n_c} \mapsto \mathbb R^{n}$ maps a state $\subs_i\in \mathbb R^{n_i+n_c}$ to a set of states $\state\in \mathbb R^{n}$:
\begin{align*} 
\proji_{\state,i}(\subs_i) :=\{\state \in \tset: (\state_i,\state_c)=\subs_i\}.
\end{align*}
A \textbf{control back projection operator} $\proji_{\ctrl,i}: \mathbb R^{m_i+m_c} \mapsto \mathbb R^{m}$ maps a control $v_i\in \mathbb R^{m_i+m_c}$ to a set of controls $\ctrl\in \mathbb R^{m}$:
\begin{align*} 
\proji_{\ctrl,i}(v_i) :=\{\ctrl \in \cset: (\ctrl_i,\ctrl_c)=v_i\}.
\end{align*}
\end{definition}
Following~\cite{zheng2023decomp}, we define the ACSS for the TV-CLVF.
\begin{definition}
\label{def:ACSS}
The ACSS of~\eqref{eqn:tv_clvf} is the set of all the control signals such that the corresponding trajectory achieves the same value $\tclvf (\state,\tinit)$:

\begin{align*}
\ACSS(\state,\tinit) = \{& \csig \in \cfset: \\
&\max_{ \tvar \in [\tinit , \thor]} e^{\gamma( \tvar -\tinit)}  \ell ( \traj (\tvar ;\tinit,\state,\ctrl (\cdot)) ) = \tclvf (\state,\tinit) \}
\end{align*}
\end{definition}

For numerical computation, define the ACS given a time step $\tstep$
\begin{align} \label{eqn: ACS}
    \ACS(\state,\tinit) = &\{ \ctrl \in \cset: \tclvf (\state,\tinit-\tstep)  - \tclvf (\state,\tinit)  - \notag \\
   &\hspace{-2em} \bigl(  D_\state \tclvf \cdot (\dyn (\state)+ \cdyn(\state) \cdot \ctrl )+ \gamma  \tclvf (\state,\tinit) \bigl) \tstep \geq 0 \}.
\end{align}

One admissible control signal $\ctrl_a(\cdot)$ for $\state$ can be computed by concatenating the ACSs over time: $\ctrl_a(\cdot) = [\ACS (\state,\tinit), \ACS (\xi(\tvar_1),\tvar_1),...,\ACS(\xi(\thor),\thor)]$. Both the ACSS and ACS are guaranteed to be non-empty~\cite{zheng2023decomp}, and they can be similarly defined for subsystems. 



\section{Problem Formulation and Approach}
In this section, we begin with an idealized case, where the system can be decomposed into two subsystems with no shared control, and show the reconstructed value function is the CLVF for the original system. This is called \textit{exact reconstruction}. We then provide a sufficient condition of exact reconstruction for cases where subsystems are coupled through shared controls. 
Further, for systems where an exact reconstruction is impossible, we apply the concept of ACS, and show how the reconstructed value function is a Lipschitz continuous CLF. 
All theorems follow naturally to the cases where the original system is decomposed into more than two subsystems. 

Denote the TV-CLVF of the two subsystems as
\begin{align*}
    &V_{\gamma,i} (\subs_i,\tinit) = \min _{v_i \in \mathbb V_i}  \max_{ \tvar \in [\tinit , \thor]} e^{\gamma( \tvar -\tinit)}  \ell_i ( \subt_i (\tvar) ) , \quad i = 1,2, 
\end{align*}
and the CLVF for subsystems as $V^\infty_{\gamma,i} : \mathcal D_{\gamma,i} \subseteq \tsset_i \mapsto \mathbb R$,
\begin{align*}
    V^\infty_{\gamma,i}(\subs_i) = \lim_{\tinit \rightarrow -\infty } V_{\gamma,i} (\subs_i,\tinit), \quad i = 1,2. 
\end{align*}
For clarity, the CLVF computed directly for the high-dimensional system~\eqref{eqn:dyn_sys} is called the \textit{original CLVF}, while the CLVF (TV-CLVF) reconstructed from subsystems is called the \textit{reconstructed CLVF (TV-CLVF)}, denoted as $\rclvf (\state)$ ($\rtclvf (\state,\tinit)$). Using~\eqref{eqn: ACS}, the ACS can be solved:  
\begin{align} \label{eqn: ACS_halfspace}
     D_\state \tclvf(\state,\tinit)  \cdyn(\state) \cdot \ctrl &\leq  \frac{\tclvf (\state,\tinit-\tstep)  -  \tclvf (\state,\tinit)}{\tstep}  \notag \\
     & \hspace{1em} -\gamma \tclvf(\state,\tinit) - D_\state \tclvf(\state,\tinit) \dyn (\state) , 
\end{align}
which is a linear inequality for single-input systems, and a half-space for multi-input systems. Though we cannot get an analytic solution of this half-space, it is numerically easy to obtain one element of it. 

\subsection{Exact Reconstruction of the CLVF}
In this section, we assume $\ell(\state) = \max (\ell_1 (\subs_1),\ell_2 (\subs_2))$, meaning the loss function for the original system equals the reconstructed loss function from subsystems. This is possible, because~\cite{gong2024robust} relaxes the choice of the loss to be any vector norm. One possible choice is $\ell_i(\subs_i) = ||\subs_i||_\infty$. 
\begin{lemma} \label{lemma: exact_decomposition}
    Assume there are no common controls after decomposition. Then, $\rclvf (\state) = \max (V^\infty_{\gamma,1}  (\subs_1) ,V^\infty_{\gamma,2}  (\subs_2)) = \clvf (\state)$ , where $\subs_1 =  \proj_{\state,1} (\state)$,  $\subs_2 =  \proj_{\state,2} (\state)$, and  $\mathcal D_\gamma =   \proji_{\state,1}(\mathcal D_{\gamma,1}) \bigcap \proji_{\state,2}(\mathcal D_{\gamma,2}).  $
    
    \begin{proof}
Define $\rtclvf (\state, \tinit) = \max ( V_{\gamma,1} (\subs_1,\tinit),V_{\gamma,2}  (\subs_2,\tinit) )$, we first show $\rtclvf(\state, \tinit) = \tclvf(\state,\tinit)$ by contradiction. 

Assume $\rtclvf (\state, \tinit) < \tclvf(\state,\tinit) $, then there exists optimal control signals $\ctrl^*(\cdot) \in \cfset$, $v_1^*(\cdot) \in \mathbb V_1$ and  $v_2^*(\cdot) \in \mathbb V_2$, s.t. 
\begin{align} \label{lemma1:eqn1}
    &\max_{ \tvar \in [\tinit , \thor]} e^{\gamma( \tvar -\tinit)}  \ell ( \traj (\tvar;\tinit, \state,\ctrl^*(\cdot) ) ) \notag \\
    >& \max \bigl(  \max_{ \tvar \in [\tinit , \thor]} e^{\gamma( \tvar -\tinit)}  \ell_1 ( \subt_1 (\tvar;\tinit, \subs_1,v_1^*(\cdot) ) ) ,\notag \\
    & \hspace{4em} \max_{ \tvar \in [\tinit , \thor]} e^{\gamma( \tvar -\tinit)}  \ell_2 ( \subt_2 (\tvar;\tinit, \subs_2,v_2^*(\cdot) ) )  \bigl) \notag\\
    = &  \max_{ \tvar \in [\tinit , \thor]} \max \bigl(  e^{\gamma( \tvar -\tinit)}  \ell_1 ( \subt_1 (\tvar;\tinit, \subs_1,v_1^*(\cdot) ) ) , \notag\\
    & \hspace{6em}  e^{\gamma( \tvar -\tinit)}  \ell_2 ( \subt_2 (\tvar;\tinit, \subs_2,v_2^*(\cdot) ) )  \bigl)\notag \\
    & =  \max_{ \tvar \in [\tinit , \thor]} e^{\gamma( \tvar -\tinit)} \max \bigl(   \ell_1 ( \subt_1 (\tvar;\tinit, \subs_1,v_1^*(\cdot) ) ) ,\notag \\
    & \hspace{9em}   \ell_2 ( \subt_2 (\tvar;\tinit, \subs_2,v_2^*(\cdot) ) )  \bigl).
\end{align}
Since the two subsystems have no shared controls, we could reconstruct one trajectory for~\eqref{eqn:dyn_sys} using the back projection: $\bar \traj (\tvar) = \proji_{\state,1} (\subt_1 (\tvar)) \bigcap   \proji_{\state,2}(\subt_2 (\tvar))$. Since the shared states evolve independent of the controls applied,~\eqref{lemma1:eqn1} becomes
\begin{align*}
    &\max_{ \tvar \in [\tinit , \thor]} e^{\gamma( \tvar -\tinit)}  \ell ( \traj (\tvar;\tinit, \state,\ctrl^*(\cdot) ) )  \\
     > & \max_{ \tvar \in [\tinit , \thor]} e^{\gamma( \tvar -\tinit)} \max \bigl(   \ell_1 ( \subt_1 (\tvar;\tinit, \subs_1,v_1^*(\cdot) ) ) , \\
    & \hspace{9em}   \ell_2 ( \subt_2 (\tvar;\tinit, \subs_2,v_2^*(\cdot) ) )  \bigl)\\
    = & \max_{ \tvar \in [\tinit , \thor]} e^{\gamma( \tvar -\tinit)}  \ell ( \bar \traj (\tvar;\tinit, \state,\ctrl^*(\cdot) ) ).
\end{align*}
This means $\ctrl^*(\cdot)$ is not optimal, which is a contradiction. 

Similarly, assume $\rtclvf(\state, \tinit) > \tclvf(\state,\tinit) $, then there also exists optimal control signals $\ctrl^*(\cdot) \in \cfset$, $v_1^*(\cdot) \in \mathbb V_1$ and  $v_2^*(\cdot) \in \mathbb V_2$, s.t.
\begin{align} \label{lemma1:eqn2}
    &\max_{ \tvar \in [\tinit , \thor]} e^{\gamma( \tvar -\tinit)}  \ell ( \traj (\tvar;\tinit, \state,\ctrl^*(\cdot) ) ) \notag \\
     < & \max_{ \tvar \in [\tinit , \thor]} e^{\gamma( \tvar -\tinit)} \max \bigl(   \ell_1 ( \subt_1 (\tvar;\tinit, \subs_1,v_1^*(\cdot) ) ) , \notag\\
    & \hspace{9em}   \ell_2 ( \subt_2 (\tvar;\tinit, \subs_2,v_2^*(\cdot) ) )  \bigl).
\end{align}
This is obtained following the same procedure as~\eqref{lemma1:eqn1}. However, projecting $\traj (\tvar) $ into the two subsystems spaces, we get $\bar \subt_1 (\tvar) = \proj_{\state,1} (\traj(\tvar))$,  $\bar \subt_2 (\tvar) = \proj_{\state,2} (\traj(\tvar))$, and we have 
\begin{align*}
    &\max_{ \tvar \in [\tinit , \thor]} e^{\gamma( \tvar -\tinit)}  \ell ( \traj (\tvar;\tinit, \state,\ctrl^*(\cdot) ) )\\
    = & \max_{ \tvar \in [\tinit , \thor]} e^{\gamma( \tvar -\tinit)} \max \bigl(   \ell_1 ( \bar \subt_1 (\tvar )),  \ell_2 ( \bar \subt_2 (\tvar) ) \bigl).
\end{align*}
Using~\eqref{lemma1:eqn2}, we have 
\begin{align*}
     &  \max_{ \tvar \in [\tinit , \thor]} e^{\gamma( \tvar -\tinit)} \max \bigl(   \ell_1 ( \bar \subt_1 (\tvar )),  \ell_2 ( \bar \subt_2 (\tvar) ) \bigl)\\
     < & \max_{ \tvar \in [\tinit , \thor]} e^{\gamma( \tvar -\tinit)} \max \bigl(   \ell_1 ( \subt_1 (\tvar;\tinit, \subs_1,v_1^*(\cdot) ) ) , \notag\\
    & \hspace{9em}   \ell_2 ( \subt_2 (\tvar;\tinit, \subs_2,v_2^*(\cdot) ) )  \bigl),
\end{align*}
which shows that either $v_1^*$ or $v_2^*$ is not optimal. Combined, we have shown that $\forall t \leq 0$, $\state \in \mathbb R^n$, $\tclvf(\state,\tinit) = \rtclvf(\state,\tinit)$. Next we show that $\bar {V}_\gamma^\infty(\state) = \clvf(\state)$.

Since the CLVF exists for two subsystems, we have
\begin{align*}
    \clvf(\state) & = \lim_{\tinit \rightarrow -\infty} \tclvf(\state,\tinit)\\
    & =\lim_{\tinit \rightarrow -\infty} \rtclvf(\state, \tinit) \\
    & = \lim_{\tinit \rightarrow -\infty} \max \bigl( V_{\gamma,1} (\subs_1,\tinit),V_{\gamma,2}  (\subs_2,\tinit) \bigl) \\
    & =  \max \bigl( \lim_{\tinit \rightarrow -\infty} V_{\gamma,1} (\subs_1,\tinit), \lim_{\tinit \rightarrow -\infty} V_{\gamma,2}  (\subs_2,\tinit) \bigl)  \\
    & = \max \bigl(  V^\infty_{\gamma,1} (\subs_1),  V^\infty_{\gamma,2}  (\subs_2) \bigl)= \tclvf(\state).
\end{align*}
This also shows that $\mathcal D_\gamma \supseteq \proji_{\state,2}(\mathcal D_{\gamma,1})  \bigcap  \proji_{\state,2}(\mathcal D_{\gamma,2})  $. 

Assume $\exists \state \notin \mathcal D_\gamma$, $\state \in \proji_{\state,2}(\mathcal D_{\gamma,1})  \bigcap  \proji_{\state,2}(\mathcal D_{\gamma,2})$. This means the optimal trajectory $\traj$ starting from $\state$ satisfies the following two condition simultaneously: 1) $\clvf(\state) $ is infinite and 2) both $ V^\infty_{\gamma,1}  (\subs_1)$ and $ V^\infty_{\gamma,2}  (\subs_2)$ are finite. This is a contradiction. Therefore, we have $\mathcal D_\gamma = \proji_{\state,2}(\mathcal D_{\gamma,1})  \bigcap  \proji_{\state,2}(\mathcal D_{\gamma,2})  $.

\end{proof}
\end{lemma}

Lemma~\ref{lemma: exact_decomposition} shows that if there exists no shared control, the reconstructed value function is the CLVF, and their domains are also the same. This means system~\eqref{eqn:dyn_sys} can be exponentially stabilized to the origin from $\mathcal D_\gamma$. 

\begin{remark}
    Numerically, Lemma~\ref{lemma: exact_decomposition} is easy to implement: after decomposition, it is easy to check if all subsystems have shared controls. If Lemma~\ref{lemma: exact_decomposition} is applicable, we simply compute the CLVF for the subsystems and take the max to reconstruct the CLVF for the original system. 
\end{remark}

The essence of this proof is that the optimal control signal for the original system can always be reconstructed from the subsystems' optimal control signals, i.e. $\ctrl^*(\cdot) = [v_1^*(\cdot);v_2^*(\cdot)]$, and that the projection and back projection between the original system's trajectory and subsystems' trajectories are both unique. This is guaranteed because there is no shared control. However, many practical systems cannot be decomposed into fully independent subsystems. From Remark~\ref{remark:shared_state}, if there are shared controls, there must also be shared states. If two subsystems require conflicting values of the shared control(s), the back projection of subsystems' trajectories cannot reconstruct the original system's trajectory. 

For the subsystems, denote the ACSS as $\ACSSo (\subs_1,\tinit)$ and $\ACSSt (\subs_2,\tinit)$, and their shared control component as $\ACSSo^c (\subs_1,\tinit)$ and $\ACSSt^c (\subs_2,\tinit)$. The ACSS for~\eqref{eqn:dyn_sys} can be reconstructed as 
\begin{align*}
    &\rACSS(\state,\tinit) = \proji_{\ctrl,1}(\ACSSo(\proj_{\state,1} (\state),\tinit) ) \bigcap \proji_{\ctrl,2} ( \ACSSt(\proj_{\state,2} (\state),\tinit) )\\
    &\rACSS^{c}(\state,\tinit) = \ACSSo^c(\proj_{\state,1} (\state),\tinit) \bigcap \ACSSt^c(\proj_{\state,2} (\state),\tinit).
\end{align*}
Note that though $\ACSSo^c (\subs_1,\tinit)$ and $\ACSSt^c (\subs_2,\tinit)$ are non-empty, $ \ACSS^{c}(\state,\tinit) $ might be empty for some $\state$. 
\begin{remark}
    One can see that $\rACSS(\state,\tinit)$ is non-empty if and only if $\rACSS^c(\state,\tinit)$ is non-empty.
\end{remark}

Now, we provide one sufficient condition on exact reconstruction for the case with shared controls.  
\begin{theorem} \label{thrm: exact_decomposition}
     Let $\rclvf (\state) = \max (V^\infty_{\gamma,1}  (\subs_1) ,V^\infty_{\gamma,2}  (\subs_2))$ where $\subs_1 =  \proj_{\state,1} (\state)$,  $\subs_2 =  \proj_{\state,2} (\state)$. Assume  $\rACSS^{c}(\state,\tinit)$ is non-empty for all $t \leq 0$ and $ \state \in \mathcal S_{\gamma}$, where $\mathcal S_{\gamma}$ is some level set of $\rclvf (\state)$. Then, $\rclvf (\state)  = \clvf (\state)$ on $\mathcal S_\gamma $.
     
\begin{proof}
The proof can be obtained similarly to Lemma~\ref{lemma: exact_decomposition}. If $\rACSS^c (\state,\tinit)$ is non-empty for all $\tinit\leq 0 $, $\state \in \mathcal S_\gamma$, then given the subsystems' trajectories $\subt_1 (\tvar)$ and $\subt_2 (\tvar)$, a trajectory of~\eqref{eqn:dyn_sys} can again be reconstructed as 
\begin{align*} 
    \bar \traj (\tvar) =\proji_{\state,1}(\subt_1 (\tvar)) \bigcap  \proji_{\state,2} (\subt_2 (\tvar)).
\end{align*}
Similarly, given a trajectory $\traj (\tvar)$ of~\eqref{eqn:dyn_sys}, the two subsystems' trajectories can be obtained as
\begin{align*} 
\bar \subt_1 (\tvar) = \proj_{\state,1}(\traj(\tvar)), \quad \bar \subt_2 (\tvar) = \proj_{\state,1}(\traj(\tvar)).
\end{align*}
Therefore, assume $\rtclvf (\state, \tinit) < \tclvf(\state,\tinit)$, we can obtain~\eqref{lemma1:eqn1}, and assume $\rtclvf (\state, \tinit) > \tclvf(\state,\tinit)$, we can obtain~\eqref{lemma1:eqn2}. In other words, $\forall \state \in \mathcal S_\gamma$, we again construct two contradictions to show that $\rtclvf(\state, \tinit) < \tclvf(\state,\tinit) $ and $\rtclvf(\state, \tinit) > \tclvf(\state,\tinit)$ cannot happen. The remaining steps are identical to the proof of Lemma~\ref{lemma: exact_decomposition}. 
\end{proof}
\end{theorem}

A direct implication of Theorem~\ref{thrm: exact_decomposition} is that the reconstructed ACSS (if not empty) is always equal to the ACSS directly computed for~\eqref{eqn:dyn_sys}. In other words, if $\rACSS^{c}(\state,\tinit)$ is non-empty for all $t \leq 0$ and $ \state \in \mathcal S_{\gamma}$, then $\rACSS^c(\state,\tinit) = \ACSS^c(\state,\tinit)$. We provide the following standard results without proof.
\begin{proposition} \label{prop:exp_stable}
    System~\eqref{eqn:dyn_sys} can be exponentially stabilized to the origin from $ {\mathcal S}_\gamma$.

\end{proposition}
\begin{algorithm} [t]
\caption{Exact reconstruction with shared controls}\label{algo:exact_decomposition}
\begin{algorithmic}[1]
\Require  Subsystems TV-CLVF $V_{\gamma,i}$, convergence time $T$, time step $\tstep$ \\

\textbf{Initialization:} $\tinit \gets T$ \\
\textbf{Output:} $\clvf (\state)$ and domain $\mathcal S_\gamma$ 
\State Reverse time of TV-CLVFs
\While{$\tinit \neq \thor$} 

\State Compute the subsystems ACS $\mathcal{U}_{\text{a},1}$, $\mathcal{U}_{\text{a},2}$ using~\eqref{eqn: ACS_halfspace}
\State Get the \textbf{shared control component} $\mathcal{U}^c_{\text{a},1}, \mathcal{U}^c_{\text{a},2}$ and other component $\mathcal{U}^o_{\text{a},1},\mathcal{U}^o_{\text{a},2}$
\State$ \ACS^{c}(\state,\tinit) \gets \ACSo^c(\text{proj}_1 (\state),\tinit) \bigcap \ACSt^c(\text{proj}_2 (\state),\tinit)$
\If{$\ACS^{c}(\state,\tinit)$ empty} Break 
\ElsIf{$\ACS^{c}(\state,\tinit)$ not empty}
\State $u_c \in \ACS^{c}(\state,\tinit)$, $u_1 \in \mathcal{U}^o_{\text{a},1}$,$u_2 \in \mathcal{U}^o_{\text{a},2}$
\State $u \gets [u_1,u_2,u_c]$
\State $\state \gets \state+ [\dyn(\state) + \cdyn(\state) \cdot u] \tstep $ 
\State $\tinit \gets \tinit + \tstep$
\EndIf
\EndWhile
\State $\clvf(\state) \gets \max (V^\infty_{\gamma,1}  (\subs_1,\thor) ,V^\infty_{\gamma,2} (\subs_2,\thor) )$, add $\state$ to $\mathcal T_\gamma$
\State $\mathcal S_\gamma \gets$ largest level set of $\clvf(\state)$ contained in $\mathcal T_\gamma$
\end{algorithmic}
\end{algorithm}

The numerical implementation of Theorem~\ref{thrm: exact_decomposition} is shown in Alg.~\ref{algo:exact_decomposition}. The convergence time for the subsystems may differ, therefore, after one TV-CLVF converges, we repeat its final value until the other subsystem's CLVF converges. Further, since the TV-CLVF is computed backward in time, we reverse the time vector of the TV-CLVF ($\tclvf(\state,\tinit) = \tclvf(\state, T - \tinit)$ when used to forward propagate the trajectory. 

Alg.~\ref{algo:exact_decomposition} has a maximum iteration $\frac{T}{\tstep}$. It first computes the ACS from the subsystem's value function and reconstructs the ACS for~\eqref{eqn:dyn_sys}. Then, it checks if the shared control component of the ACS is empty. If not empty, it applies one random control input from ACS, updates the state, and repeats. For any initial state whose ACS is non-empty along the whole trajectory, concatenating the control used, we get one admissible control signal. This state is added to a set $\mathcal T_\gamma$, and $\mathcal S_\gamma$ is the largest level set of $\rclvf(\state)$ contained in $\mathcal T_\gamma$.  Note that we take the maximum between $V^\infty_{\gamma,1}  (\subs_1,\thor) $ and $V^\infty_{\gamma,2} (\subs_2,\thor)$. This is because we reverse the time vector at line 3. 


\subsection{Inexact Reconstruction: ACS and lipschitz continuous CLF}

Lemma~\ref{lemma: exact_decomposition} and Theorem~\ref{thrm: exact_decomposition} both provide exact reconstructions of the CLVF, but on different domains. In this section, we generalize to the case in which the subsystems have shared control, and therefore the sufficient condition in Theorem~\ref{thrm: exact_decomposition} is not met. 
Though exact reconstruction of a CLVF in the original state space is impossible, we can still reconstruct a Lipschitz continuous CLF and it ROES using the subsystems' CLVFs and ACSs. 
We remove the assumption that $\ell(\state) = \max (\ell_1 (\subs_1),\ell_2 (\subs_2))$. 

A key difference compared to the exact reconstruction is that we focus on the $\ACS (\state)$ for $\clvf (\state)$, instead of $\ACSS (\state,\tinit)$ for $\tclvf (\state,\tinit)$. That is to say, we do not care about how the subsystems TV-CLVF evolve, but only care about its CLVF. Note that for the CLVF, its ACS becomes time-independent: 
\begin{align} \label{eqn: ACS_clvf}
    \ACS(\state) = \{ \ctrl \in \cset: D_\state \clvf[\dyn (\state)+ \cdyn(\state) \cdot \ctrl]  \leq -\gamma \clvf(\state) \}.
\end{align}
Denote the ACSs for two subsystems as $\ACSo(\subs_1)$, $\ACSt(\subs_2)$, and their shared control components as $\ACSo^c (\subs_1)$ and $\ACSt^c (\subs_2)$. The ACS for~\eqref{eqn:dyn_sys} can be reconstructed as
\begin{align*}
    &\rACS(\state) = \proji_{\ctrl,1}(\ACSo(\proj_{\state,1} (\state)) ) \bigcap \proji_{\ctrl,2} ( \ACSt(\proj_{\state,2} (\state)) )\\
    &\rACS^{c}(\state) = \ACSo^c(\proj_{\state,1} (\state)) \bigcap \ACSt^c(\proj_{\state,2} (\state)).
\end{align*}
\begin{theorem}
     \label{thrm: inexact_decomposition}
     Let $\rclvf (\state) = V^\infty_{\gamma,1}  (\subs_1) + V^\infty_{\gamma,2}  (\subs_2)$ where $\subs_1 =  \proj_{\state,1} (\state)$,  $\subs_2 =  \proj_{\state,2} (\state)$. Assume $\rACS^{c}(\state)$ is non-empty for all $ \state \in \bar {\mathcal S}_{\gamma}$, where $\bar {\mathcal S}_{\gamma}$ is some level set of $\rclvf$. Then, $\rclvf$ is a Lipschitz continuous local CLF for~\eqref{eqn:dyn_sys} on $\bar {\mathcal S}_\gamma $.
     
     \begin{proof}
From~\eqref{eqn:Vdot}, there exist $v_1 \in \ACSo(\subs_1) $ and $v_2 \in \ACSt(\subs_2) $, such that 
\begin{align*}
    &D_{\state} \rclvf (\state) \cdot (\dyn (\state)+ \cdyn(\state) \cdot \ctrl)\\
    =&  D_{\state} V^\infty_{\gamma,1}  (\subs_1) \cdot ( \dyn_1(\subs_1)+g_1(\subs_1)v_1) + \\
    & \hspace{1em} D_{\state} V^\infty_{\gamma,2}  (\subs_2) \cdot ( \dyn_2(\subs_2)+g_2(\subs_2)v_2) \\
     \leq &-\gamma V^\infty_{\gamma,1}  (\subs_1)  -\gamma V^\infty_{\gamma,2}  (\subs_2) \\
     =& -\gamma \rclvf (\state),
\end{align*}
where $\ctrl = \proji_{\ctrl,1}(v_1) \bigcap \proji_{\ctrl,2} (v_2 )$. In other words, $\forall \state \in \bar {\mathcal S}_\gamma$, there exists some $\ctrl \in \cfset$ s.t. the Lie derivative along~\eqref{eqn:dyn_sys} is smaller than $-\gamma \rclvf (\state)$. 

Further, since both $V^\infty_{\gamma,1}$ and $V^\infty_{\gamma,2}$ are positive definite and all level sets are closed, $\rclvf (\state) = V^\infty_{\gamma,1}+V^\infty_{\gamma,2}$ must also be positive definite and have closed level sets. We conclude that $\rclvf (\state)$ is a local CLF. 
\end{proof}
\end{theorem}

Similar to Proposition~\ref{prop:exp_stable}, system~\eqref{eqn:dyn_sys} can be exponentially stabilized to the origin from $\bar {\mathcal S}_\gamma$.

\begin{remark}
    We use the summation instead of the maximum to reconstruct the value function. This is because using summation does not introduce additional non-differentiable points, whereas for the maximum, the reconstructed value function is non-differentiable for all $\state $ where $V^\infty_{\gamma,1}  (\subs_1) = V^\infty_{\gamma,2}  (\subs_2)$. 
\end{remark}

\subsection{Optimal QP Controller }

For both exact and inexact reconstruction, the controller that guarantees exponential stabilizability can be synthesized by solving the following QP, 
\begin{align} \label{eqn:QP}
   & \hspace{5em} \min_{\ctrl \in \cset} ||\ctrl - \ctrl_r|| \notag \\
    &\text{ s.t. } D_\state \rclvf[\dyn (\state)+ \cdyn(\state) \cdot \ctrl]  \leq -\gamma \rclvf (\state) .
\end{align}
This QP is guaranteed to be feasible on $\mathcal D_\gamma$ (or $ {\mathcal S}_\gamma $) \cite{gong2022constructing}.
For systems whose CLVFs can be exactly reconstructed, the constraints in~\eqref{eqn:QP} can either be solved directly from the reconstructed CLVF, or from subsystems CLVFs. 
For systems in which the CLVFs cannot be exactly reconstructed, the control can be determined by solving~\eqref{eqn:QP} with the reconstructed CLF. 



\section{Numerical Example}
In this section, we provide numerical examples that validate our theory. This includes examples in 2D, 3D, and 10D. For simplicity, in all examples, the loss function is chosen to be the infinity norm. All simulations are conducted using MATLAB and tool boxes \cite{ian2005levelset,chenoptimal}.

\begin{figure}
    \centering
    \includegraphics[width=\columnwidth]{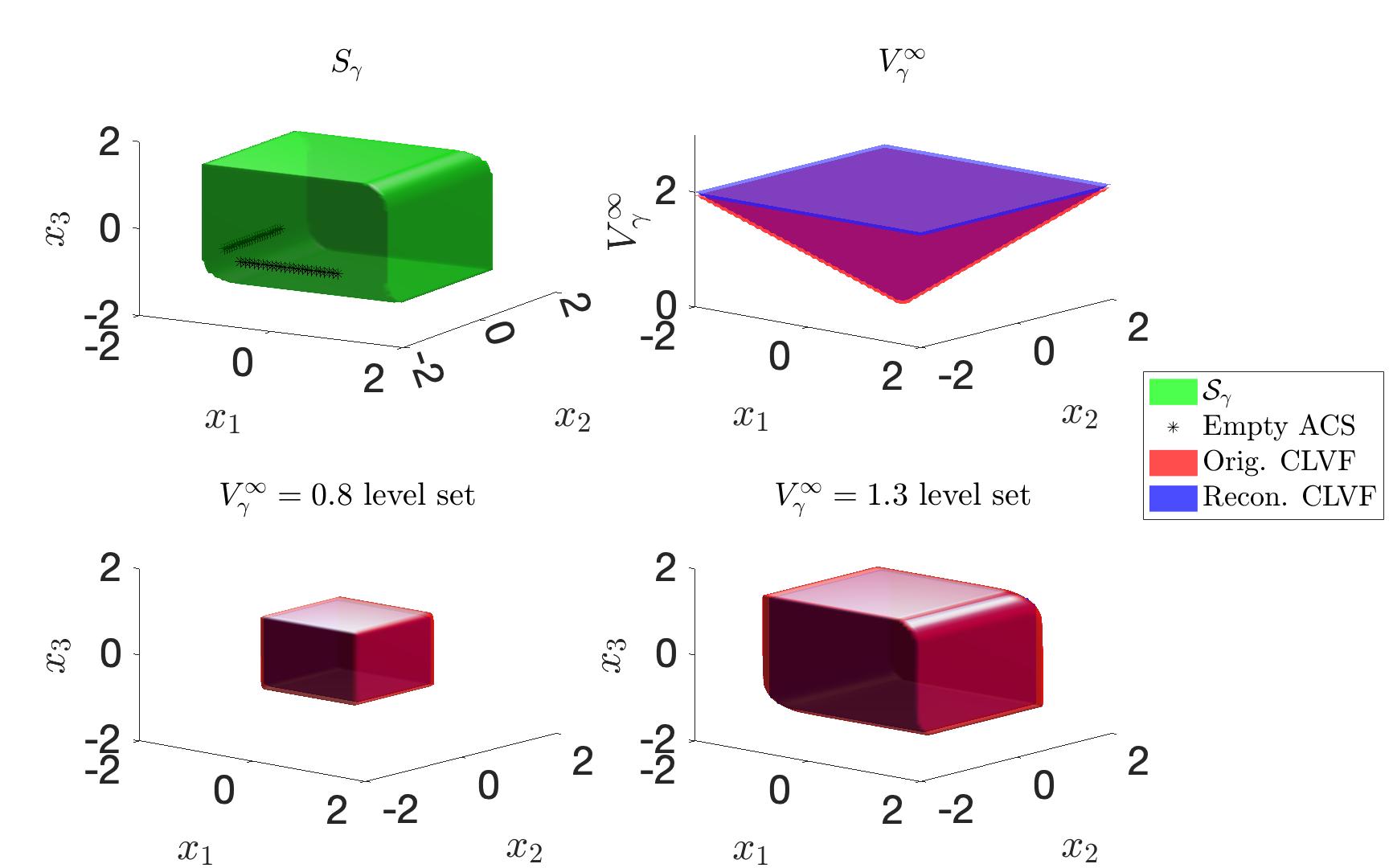}
    \caption{Original and reconstructed CLVFs with $\gamma=0.1$ for system~\eqref{eqn: multlinear3D}. Top left, $S_\gamma$ for the reconstructed ACS produced by Alg.~\ref{algo:exact_decomposition}. Top right, comparison of the reconstructed and original CLVF projected into $x_1-x_2$ plane. Bottom, different level sets of original and reconstructed CLVFs. The reconstructed CLVF is identical to the original CLVF in $\mathcal S_\gamma$, validating Theorem~\ref{thrm: exact_decomposition}. The computation time for the original system's CLVF is 102s on a grid of $[-2,-2,-2]$ to $[2,2,2]$ with 61 grid points on each dimension. On the same grid, the computation time for the subsystems' CLVFs is 3.37s, and the computation of $\mathcal S_\gamma$ takes 0.86s. Combined, the decomposition speeds computation time by 30x.}
    \label{fig:multlinear3D}
\end{figure}

\subsection{2D System (Lemma~\ref{lemma: exact_decomposition})}
Consider the following 2D system 
\begin{align} \label{eqn: nonlinear2D}
    \dot{x}_1 = x_1^2+u_1, \quad \dot{x}_2 = x_2+u_2,   
\end{align}
where $u_1 \in [-4,4]$ and $u_2 \in [-1,1]$. The two states evolve purely on their own, therefore each state's dynamics is one subsystem. We compute the CLVF for both subsystems and also directly for the original system with $\gamma = 0.1$ and $0.3$. The results are shown in Fig.~\ref{fig:2Dsys}.

\subsection{3D System (Theorem~\ref{thrm: exact_decomposition})}

Consider the following 3D system 
\begin{align} \label{eqn: multlinear3D}
     \dot{x_1} = x_3 + u_1, \hspace{2em} \dot{x_2} =  x_3 + u_2, \hspace{2em} \dot{x_3} = u_3,
\end{align}
where $u_1, u_2\in[-1,1]$ and $u_3 \in [-0.5,0.5]$ are the control inputs. 
Take $z_1 = [x_1,x_3]$, $z_2 = [x_2,x_3]$, and $v_1 = [u_1,u_3]$, $v_2 = [u_2,u_3]$, we can decompose the system into two 2D subsystems. We can apply Alg.~\ref{algo:exact_decomposition} to find the domain $S_\gamma$, such that if initialized inside, the system will stabilize to its equilibrium point. The results are shown in Fig.~\ref{fig:multlinear3D}. 


\begin{figure}
    \centering
    \includegraphics[width=\columnwidth]{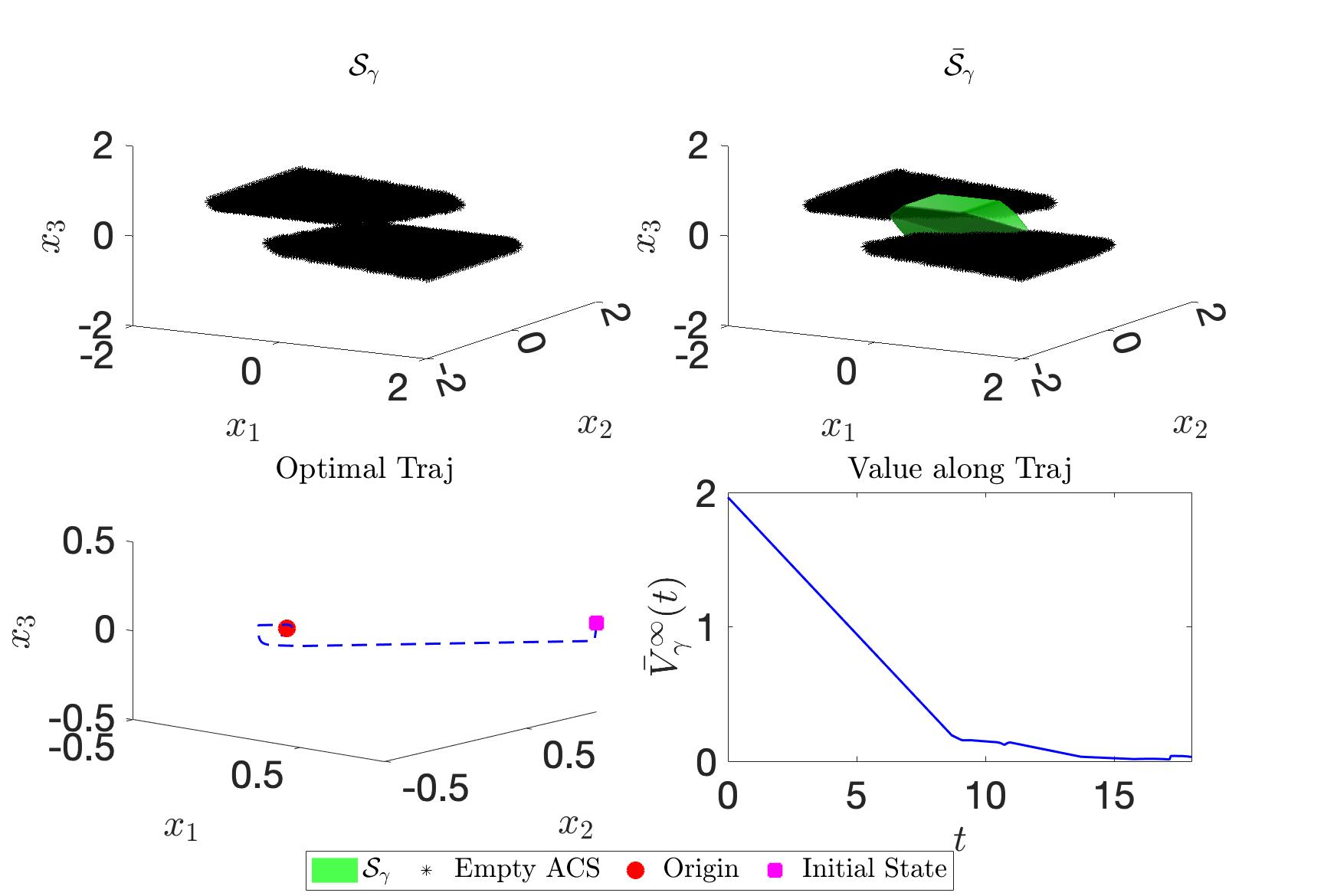}
    \caption{ Top: The black stars denote the states with empty ACS (and ACSS). Left, $\mathcal S_\gamma $ (consists only of the origin) computed from Alg.~\ref{algo:exact_decomposition}, shown in green. Right, $\bar {\mathcal S}_\gamma $. Bottom left: trajectory simulation using the QP~\eqref{eqn:QP}. Bottom right: decay of the CLF value along the trajectory. When initialized inside $\bar {\mathcal S}_\gamma$, the system can be stabilized to the origin, validating Theorem ~\ref{thrm: inexact_decomposition}. The computation time for the subsystems' CLVFs is 15.5s with a grid from $[-2,-2]$ to $[2,2]$ and 101 nodes on each dimension. The computation of $\mathcal S_\gamma$ takes 3.26s.}
    \label{fig:linear3D}
\end{figure}

\subsection{3D System (Theorem~\ref{thrm: inexact_decomposition})}
Consider the following 3D system 
\begin{align} \label{eqn: linear3D}
    \dot{x_1} = x_3 , \hspace{2em} \dot{x_2} =  x_3, \hspace{2em} \dot{x_3} = u,
\end{align}
where $u_1\in[-1,1]$ is the control inputs. Take $z_1 = [x_1,x_3]$, $z_2 = [x_2,x_3]$, and $v_1 = [u_1,u_3]$, $v_2 = [u_2,u_3]$, we can decompose the system into two 2D subsystems. Although reconstruction to the exact CLVF is not possible, we can still validate the reconstructed lipschitz continuous CLF through a trajectory simulation inside $S_\gamma$. The results are shown in Fig.~\ref{fig:linear3D}. 


\subsection{10D Quadrotor}
Consider the 10D quadrotor system: 
\begin{align}
    &\dot{x}_1 = x_2 , \hspace{2mm} \dot{x}_2 =  g\tan{x_3}, \hspace{2mm} \dot{x_3} = -d_1 x_3 + x_4, \notag \\ 
    &\dot{x}_4 =  -d_0 x_3 + n_0 u_1, \hspace{2mm} \dot{x}_5 =  x_6 , \hspace{2mm}  \dot{x}_6 = g\tan{x_7}, \notag \\
    & \dot{x}_7 = -d_1x_7 + x_8, \hspace{2mm} \dot{x}_8 = -d_0x_7 + n_0u_2, \notag \\
    \
    & \dot{x}_9 = x_{10} , \hspace{2mm} \dot{x}_{10} = u_3, \label{eq: 10D_Quad}
\end{align}
where $(x,y,z)$ denote the position, $(x_2, x_6, x_{10})$ denote the velocity, $(x_3, x_7)$ denote the pitch and roll, $(x_4, x_8)$ denote the pitch and roll rates, and $(u_1, u_2, u_3)$ are the controls. The system parameters are set to be $d_0 = 10, d_1 = 8, n_0 = 10, k_T = 0.91, g = 9.81$, $|u_1|, |u_2| \leq \pi/4$, $u_3\in [-1, 1]$.

\begin{figure}
    \centering
    \includegraphics[width=\columnwidth]{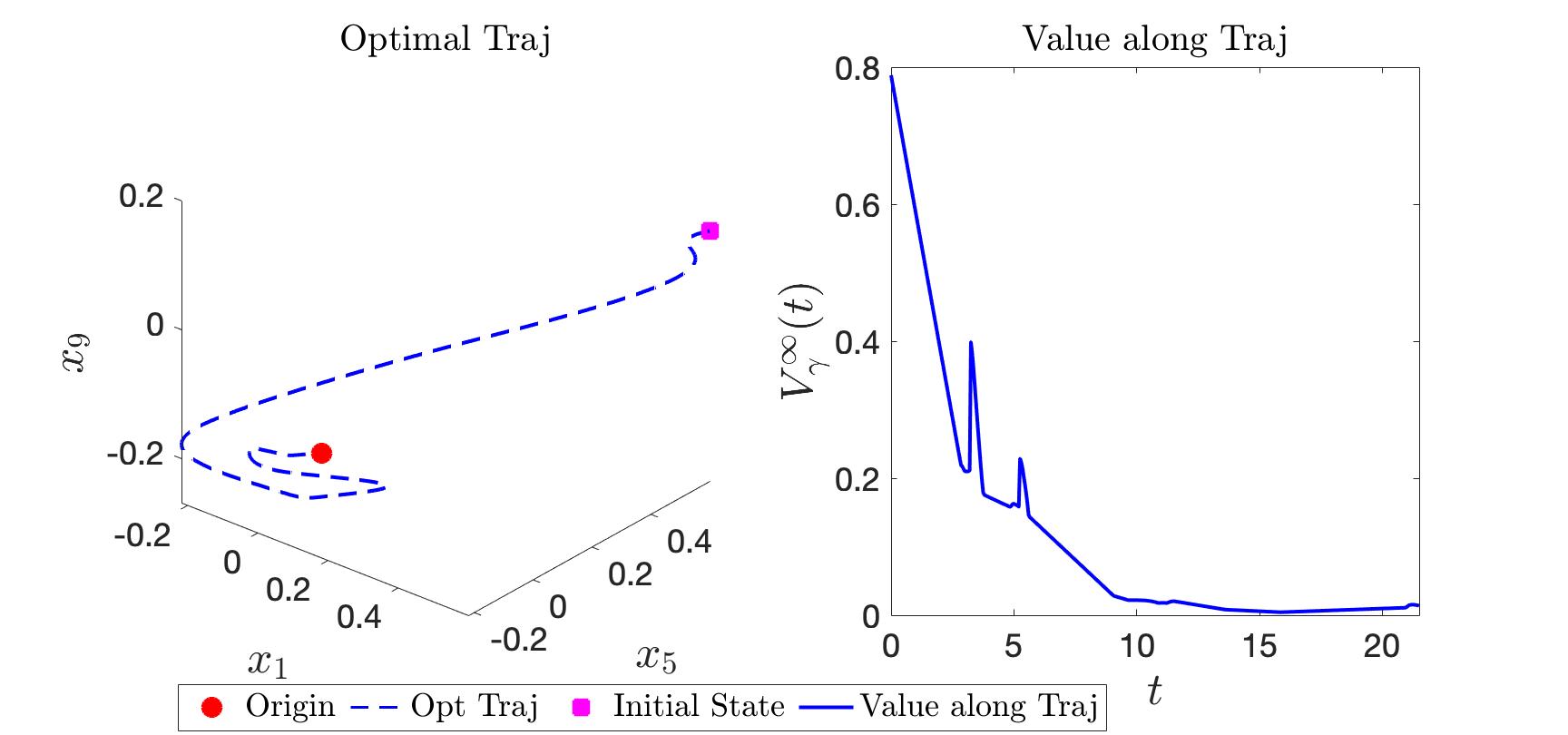}
    \caption{Since the CLVF is a 10D function, we do not visualize it, but show instead a trajectory that is stabilized to the origin and the value decay along this trajectory. The computation for the X/Y subsystem takes 1374.51s, and 56.63s for the Z subsystem. The direct computation for the 10D CLVF is not tractable.}
    \label{fig:10D}
\end{figure}

This 10D system can be decomposed into three subsystems: $X$ subsystem with states $[x_1,x_2,x_3,x_4]$,  $Y$ subsystem with states $[y,x_6,x_7,x_8]$, and $Z$ subsystem with states $[x_9,x_{10}]$. It can be verified that all three subsystems have an equilibrium point at the origin. Further, there's no shared control or states among subsystems, therefore, Lemma~\ref{lemma: exact_decomposition} can be used to exactly reconstruct the CLVF using $\bar \ell(\state) = ||\state||_\infty$. The result is shown in Fig.~\ref{fig:10D}.

\section{Conclusion}
In this paper, we presented an exact and inexact reconstruction method for CLVFs by leveraging ACS, providing a corresponding QP controller, and identifying the domain that can exponentially stabilize the system. We validate our method with a nonlinear 2D system with no shared control for exact reconstruction, as well as a multi-input linear 3D system for the shared control case. Further, we validate the inexact reconstruction method with the single-input linear 3D system by finding the corresponding domain. Lastly, a 10D quadrotor example is provided, showing numerical efficiency. 

Through this method, we can scalably compute higher-dimensional CLVFs for a class of nonlinear systems. Future work includes validating high-dimensional CLVFs through neural solvers like DeepReach~\cite{bansal2021deepreach}, applying this method to online trajectory planning problems, finding the ``smallest control invariant set'' defined in the original CLVF work~\cite{gong2022constructing}, and exploring more on the region where the reconstructed ACS is empty.

\bibliographystyle{IEEEtran}\begin{scriptsize}
\bibliography{ref}\end{scriptsize}

\end{document}